\documentclass{amsart}
\usepackage{float,mathrsfs,amsmath,amssymb,amsthm,multirow,tabls}
\usepackage{graphicx}
\usepackage[all]{xy}

\numberwithin{equation}{section}

\theoremstyle{plain}
\newtheorem{theorem}{Theorem}

\newtheorem*{acknowledgements}{Acknowledgements}

\theoremstyle{remark}
\newtheorem{remark}[theorem]{Remark}

\DeclareMathOperator{\I}{I}
\DeclareMathOperator{\II}{II}
\DeclareMathOperator{\III}{III}
\DeclareMathOperator{\IV}{IV}

\title[Jacobian fibrations on $X_3$]{Jacobian fibrations on the
  singular $K3$ surface of discriminant 3}
\date{}
\author{Kazuki Utsumi}
\address{Department of Mathematics \endgraf Graduate School
  of Science \endgraf Hiroshima University \endgraf 1-3-1 Kagamiyama,
  Higashi-Hiroshima 739-8526 \endgraf Hiroshima University }
\email{kazu-utsumi@hiroshima-u.ac.jp}

\begin{document}

\maketitle

\begin{abstract}
  In this paper we give the Weierstrass equations and the generators
  of Mordell-Weil groups for Jacobian fibrations on the singular $K3$
  surface of discriminant 3. 
\end{abstract}

\section{Introduction}

A $K3$ surface defined over the complex number field whose Picard
number equals to maximum possible number $20$ is called a
\textit{singular $K3$ surface}. Shioda and Inose~\cite{Shio-Ino}
showed that the map a singular $K3$ surface $X$ corresponds to its
transcendental lattice $T_X$ is a bijective correspondence from the
set of singular $K3$ surfaces onto the set of equivalence classes of
positive-definite even integral lattice of rank two with respect to
$SL_2(\mathbb{Z}$). The discriminant of a singular $K3$ surface $X$ is
the determinant of the Gram matrix of the transcendental lattice
$T_X$.

In this paper we study Jacobian fibrations, i.e., elliptic fibrations
with a section, on the singular $K3$ surface $X_3$ of discriminant
$3$, which corresponds to the lattice defined by $\left(
  \begin{array}{cc}
    2 & 1\\
    1 & 2
  \end{array}
\right)$ and is uniquely determined up to isomorphism. Jacobian
fibrations on $X_3$ were classified by Nishiyama~\cite{Nishiyama}. He
classified all configurations of singular fibers of Jacobian
fibrations on $X_3$ into $6$ classes and determined their Mordell-Weil
groups. Then, we give for each fibration a Weierstrass model. More
precisely, we state our main theorem.

\begin{theorem}
  Let $X_3$ be the singular $K3$ surface of discriminant $3$. For each
  Jacobian fibration in Nishiyama's list~\cite[Table~1.1]{Nishiyama},
  an elliptic parameter $u_i$, a Weierstrass equation and the
  generators of the Mordell-Weil group are given by
  Table~\ref{classes}.
\end{theorem}

\begin{table}[h]
  {\footnotesize
  \centering
  \begin{tabular}{|c|c|c|c|c|}
    \hline
    No. & sing. fibs & MWG & $u_i$ & equation and rational points\\ \hline
    \multirow{2}{*}{1} & \multirow{2}{*}{$2 \II^* + \IV$}  & \multirow{2}{*}{$0$} & \multirow{2}{*}{$\frac{2(y_2+1)}{(y_1-1)^2}$}
    \multirow{2}{*}{} & $Y^2=X^3+u_1^5(u_1-1)^2$ \\
    & & & & $O$\\ \hline 
    \multirow{2}{*}{2} & \multirow{2}{*}{$\I_{12}^* + \I_3+3\I_1$} & \multirow{2}{*}{$\mathbb{Z} / 2 \mathbb{Z}$} &
    \multirow{2}{*}{$\frac{2t^2}{(y_2+1)(y_1^2+2y_1+2y_2-1)}$} & $Y^2=X^3-2u_2(u_2^3-2)X^2+u_2^8X$\\
    & & & & $O, \; (0,0)$\\ \hline
    \multirow{3}{*}{3} & \multirow{3}{*}{$\III^* + \I_6^*+3\I_1$} &
    \multirow{3}{*}{$\langle \frac{3}{2} \rangle \oplus \mathbb{Z}/2\mathbb{Z}$} & \multirow{3}{*}{$\frac{t}{y_1^2-1}$}
    &$Y^2= X^3 +4u_3^3X^2-4u_3^3X$\\
    & & & & 2-tor.: $O, \; (0,0)$\\
    & & & & free gen. : $(1,-1)$\\ \hline
    \multirow{3}{*}{4} & \multirow{3}{*}{$\I_{18} + 6\I_1$} & \multirow{3}{*}{$\langle \frac{3}{2} \rangle \oplus \mathbb{Z}/3\mathbb{Z}$} &
    \multirow{3}{*}{$\frac{t}{y_1+y_2}$} & $Y^2=X^3+(X-u_4^6)^2$\\
    & & & & 3-tor. : $O, \; (0, \pm u_4^6)$ \\
    & & & & free gen. : $(2u_4^3, 2u_4^3+u_4^6)$ \\ \hline
    \multirow{2}{*}{5} & \multirow{2}{*}{$3\IV^*$} & \multirow{2}{*}{$\mathbb{Z}/3\mathbb{Z}$} &
    \multirow{2}{*}{$y_1$} &  $Y^2=X^3 + (u_5^2-1)^4$ \\
    & & & & $O, \; (0, \pm (u_5^2-1)^2))$ \\ \hline 
    \multirow{2}{*}{6} & \multirow{2}{*}{$\I_3^* + \I_{12}+3\I_1$} & \multirow{2}{*}{$\mathbb{Z}/4\mathbb{Z}$} &
    \multirow{2}{*}{$t$} & $Y^2=X^3-2(u_6^3-2)X^2+u_6^6X$ \\
    & & & & $O, \; (0,0), \; (u_6^3, \pm 2u_6^6)$ \\ \hline
  \end{tabular}
}
  \vspace{.1in}
  \caption{Classification of Jacobian fibrations on $X_3$}
  \label{classes}
\end{table}

We explain about Table~\ref{classes}. The fist column shows the name
of each Jacobian fibrations following Nishiyama's notaion. The second
column shows the configuration of singular fibers. Here, for example,
by $2\II^*+\IV$ meanes that the surface has two singular fibers of
type $\II^*$ and a singular fiber of ot type $\IV$ (Kodaira's
notation~\cite{Kodaira}). The third column shows the Mordell-Weil
group (MWG) of the fibration. The fourth column shows an elliptic
parameter $u_i$ of the fibration under the singular affine modell
(\ref{model}) of $X_3$. The index $i$ is the name of the
fibration. The last column shows a Weierstrass equation and rational
points corresponding to Mordell-Weil generator of the fibration, where
$O$ is the rational point corresponding to the zero of MWG.

Recently, Braun, Kimura and Watari~\cite{Isom} showed that Nishiyama's
list also gives the classification of Jacobian fibrations on $X_3$
modulo isomorphism. Thus, our and their results answer completely a
question of Kuwata and Shioda~\cite{Kuwa-Shio}.

\section{Notation}

\begin{figure}[h]
  \centering
  \includegraphics[width=10cm]{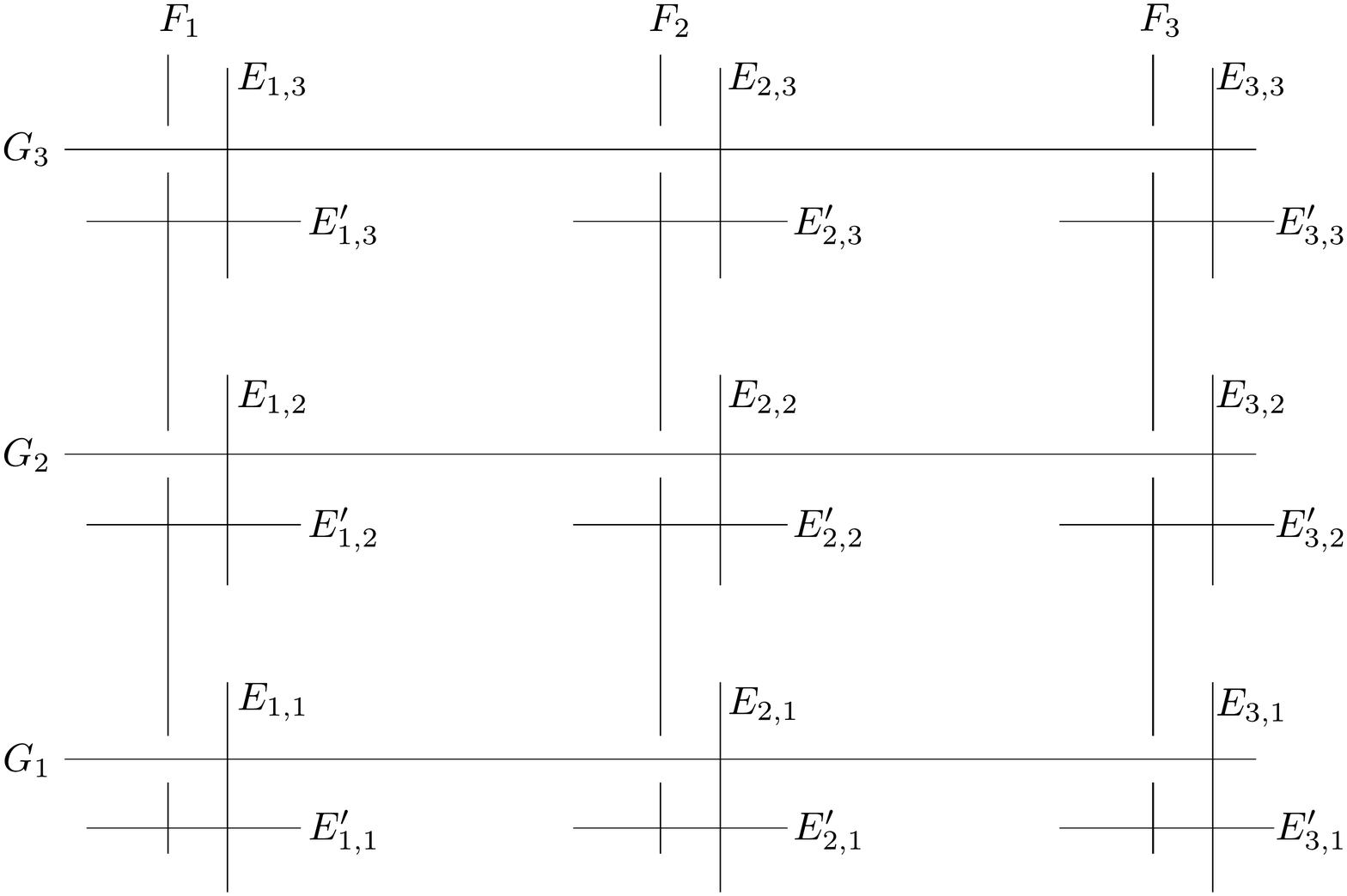}
  \caption{$(-2)$-curves}\label{config}
\end{figure}

The singular $K3$ surface $X_3$ is known as \textit{a generalized
  Kummer surface} constructed by the following. Let $C_{\omega}$ be
the complex elliptic curve with the fundamental periods $1$ and
$\omega=e^{2\pi\sqrt{-1}/3}$. Let $\sigma$ be an automorphism of
$C_{\omega} \times C_{\omega}$ defined by $\sigma(z_1, z_2) \mapsto
(\omega z_1, \omega^2 z_2)$. Then the minimal resolution of the
quotient $C_{\omega} \times C_{\omega} / \langle \sigma \rangle$ is
isomorphic to the singular $K3$ surface $X_3$ (see
\cite[Lemma~5.1]{Shio-Ino}). The automorphism $\sigma$ has the $9$
fixed points $(v_i, v_j) \, (1 \leq i, j \leq 3)$, where $\{v_i\}$ are
the fixed points of the automorphism $\sigma_1$ of $C_{\omega}$
defined by $\sigma_1(z)=\omega z$. These $9$ points $(v_i, v_j)$
correspond to the singular points $p_{ij}$ of the quotient $C_{\omega}
\times C_{\omega}/ \langle \sigma \rangle$. The minimal resolution
$X_3$ of $C_{\omega} \times C_{\omega}/ \langle \sigma \rangle$ is
obtained by replacing each $p_{ij}$ by $2$ non-singular rational
curves $E_{i,j}$ and $E_{i,j}'$ with $E_{i,j} \cdot
E_{i,j}'=1$. Moreover, $X_3$ contains $6$ non-singular rational
curves, i.e. the image $F_i$ (or $G_j$) of $\{v_i\} \times C_{\omega}$
(or $C_{\omega} \times \{v_j\}$) in $X_3$. We have the following
intersection numbers.

\begin{equation}
  \label{intsct}
  \begin{aligned}
    &F_i^2 = G_i^2 = E_{i,j}^2 =  E_{i,j}'^2 = -2, \; F_i \cdot E_{j,k} =
    G_i \cdot E_{j,k}'= F_i \cdot G_j=0,\\
    &E_{i,j} \cdot E'_{k,l}= \delta_{i,k}
    \cdot \delta_{j,l}, \quad  F_i \cdot E_{j,k}' = G_i \cdot E_{k,j} = \delta_{i,j}.
  \end{aligned}
\end{equation}
These $24$ curves on $X_3$ form the configuration of
Figure~\ref{config}.

It is well known that the elliptic curve $C_{\omega}$ has the
following Weierstrass form
\begin{equation}
  C_{\omega}: y^2=x^3+1.
\end{equation}
We denote each factor of $C_{\omega} \times C_{\omega}$ by
\begin{equation}
C_{\omega}^1 : y_1^2 = x_1^3+1, \quad C_{\omega}^2 : y_2^2 = x_2^3+1 .
\end{equation}
Then the automorphism $\sigma$ is written by
\begin{equation}
  \begin{aligned}
    \sigma \ : & \ C_{\omega}^1 \times C_{\omega}^2 \to C_{\omega}^1 \times C_{\omega}^2\\
    & (x_1, y_1, x_2, y_2) \mapsto (\omega x_1, y_1, \omega^2 x_2, y_2).
  \end{aligned}
\end{equation}
The function field $\mathbb{C}(X_3)$ is equal to the invariant
subfield of the function field $\mathbb{C}(C_{\omega}^1 \times
C_{\omega}^2)=\mathbb{C}(x_1, x_2, y_1, y_2)$ under the automorphism $\sigma$.
Then we have
\begin{equation}
\mathbb{C}(X_3) = \mathbb{C}(y_1, y_2, t), \quad t=x_1 x_2,   
\end{equation}
where $y_1, y_2,$ and $t$ are naturally regarded as functions on $X_3$ with the relation
\begin{equation}
  \label{model}
  t^3 = (y_1^2-1)(y_2^2-1).
\end{equation}
This gives a singular affine model of $X_3$. We start from the
equation to obtain a Weierstrass form for each Jacobian
fibration on $X_3$. Under the above notation, we see that the divisor of
typical functions are as follows.
\begin{equation}
  \label{func}
  \begin{aligned}
    &\big( y_1-1 \big) = 3F_2 + 2(E'_{2,1}+E'_{2,2}+E'_{2,3}) +E_{2,1}+E_{2,2}+E_{2,3}\\
    & \qquad \qquad \qquad   - \left( 3F_1 +2\left(E'_{1,1}+E'_{1,2}+E'_{1,3} \right) + E_{1,1}+E_{1,2}+E_{1,3} \right)\\
    & \\
    &\big( y_1+1 \big) = 3F_3 + 2(E'_{3,1}+E'_{3,2}+E'_{3,3}) +E_{3,1}+E_{3,2}+E_{3,3}\\
    & \qquad \qquad \qquad - \left( 3F_1 +2\left(E'_{1,1}+E'_{1,2}+E'_{1,3} \right) + E_{1,1}+E_{1,2}+E_{1,3} \right)\\
    & \\
    &\big( y_2-1 \big) = 3G_2 + 2(E_{1,2}+E_{2,2}+E_{3,2}) +E'_{1,2}+E'_{2,2}+E'_{3,2}\\
    & \qquad \qquad \qquad- \left( 3G_1 +2\left(E_{1,1}+E_{2,1}+E_{3,1} \right) + E'_{1,1}+E'_{2,1}+E'_{3,1} \right)\\ 
    & \\
    &\big( y_2+1 \big) = 3G_3 + 2(E_{1,3}+E_{2,3}+E_{3,3}) +E'_{1,3}+E'_{2,3}+E'_{3,3}\\
    & \qquad \qquad \qquad- \left( 3G_1 +2\left(E_{1,1}+E_{2,1}+E_{3,1} \right) + E'_{1,1}+E'_{2,1}+E'_{3,1} \right)\\ 
    & \\
    &\big( t \big) = F_2+E'_{2,3}+E_{2,3}+G_3+E_{3,3}+E'_{3,3}+F_3+E'_{3,2}+E_{3,2}+G_2+E_{2,2}+E'_{2,2}\\
    & \qquad \qquad - \left(E_{2,1}+E_{3,1} + 2 \left(G_1+E_{1,1}+E'_{1,1}+F_1 \right) +E'_{1,2}+E'_{1,3} \right).
  \end{aligned}
\end{equation}

\newpage

\section{Fibration 1}\label{fib1}

\begin{figure}[h]
  \centering
  \includegraphics[width=10cm]{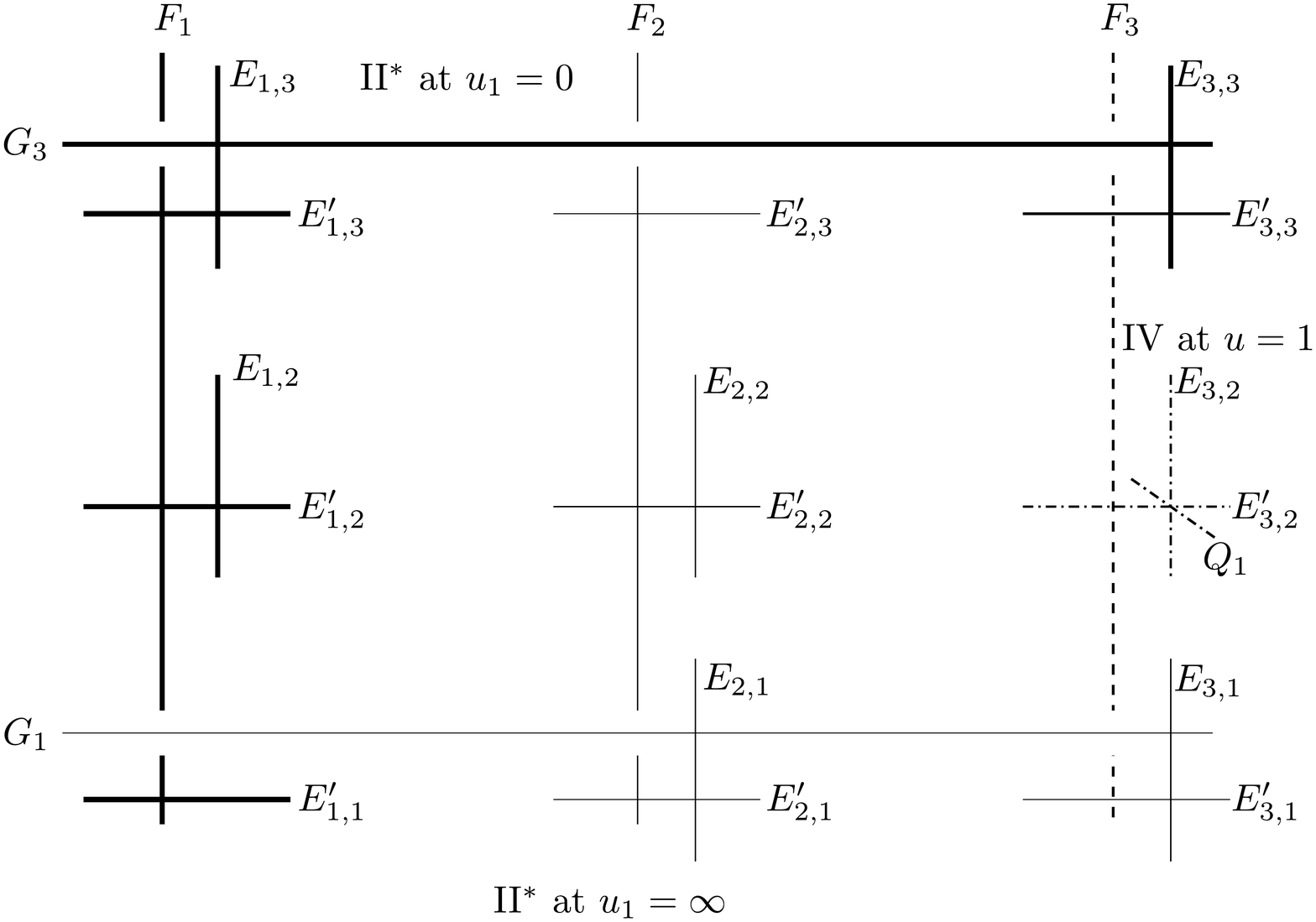}
  \caption{Fibration 1}
  \label{no1}
\end{figure}

An elliptic parameter for Fibration 1 is given by 
\begin{equation}
  \label{par1}
  u_1=\frac{2(y_1+1)}{(y_1-1)^2}.
\end{equation}
The divisor of $u_1$ is given by
\begin{equation}
  \label{div1}
  \begin{aligned}
    \left( u_1 \right) &= E'_{3,3}+2E_{3,3}+3G_3+4E_{1,3}+5E'_{1,3}+6F_1+3E'_{1,1}+4E'_{1,2}+2E_{1,2}\\
    & \quad - \left( E'_{3,1}+ 2E_{3,1}+3G_1+4E_{2,1}+5E'_{2,1}+6F_2+3E'_{2,3}+4E'_{2,2}+2E_{2,2} \right).
  \end{aligned}
\end{equation}
The zero divisor $(u_1)_0$ (the bold lines in Figure~\ref{no1}) and
the polar divisor $(u_1)_{\infty}$ (the thin lines in
Figure~\ref{no1}) are the singular fibers both of type $\II^*$. 

Eliminating the variable $y_2$ from (\ref{model}) and (\ref{par1}), we
obtain the following equation
\begin{equation}
  \label{eq1-1}
  4t^3= u_1 (y_1+1)(y_1-1)^3(u_1 y_1^2-2 u_1 y_1+ u_1 -4),
\end{equation}
which defines a plane curve over $\mathbb{C}(u_1)$ with a singularity
at $(t,y_1)=(0,1)$. Blowing up by $t=v(y_1-1)$, we have the following
equation
\begin{equation}
  \label{eq1-2}
  4v^3=u_1 (y_1+1) ( u_1 y_1^2-2u_1y_1+u_1-4),
\end{equation}
which defines a nonsingular plane cubic curve over $\mathbb{C}(u_1)$
with a rational point $(v,y_1)=(0,-1)$. Then we can convert it into a
Weierstrass form (see \cite{Jacobians} or \cite{Connell}). Since the
rational point $(v,y_1)=(0,-1)$ corresponds to the divisor $F_3$ (the
dotted line in Figure~\ref{no1}), choosing it as the zero section of
the group structure, we obtain the Weierstrass equation for Fibration
1
\begin{equation}
  \label{eq-no1}
  Y^2=X^3+u_1^5(u_1-1)^2,
\end{equation}
where the change of variables is given by
\begin{equation}
  \label{change1}
    X = \frac{\sqrt[3]{4}(u_1-1) u_1 t}{(y_1^2-1)}, \quad   Y = -\frac{u_1^2 (u_1-1) (u_1 y_1 -u_1 +2)}{y_1+1}.
\end{equation}
Besides the two singular fibers of type $\II^*$ at $u_1 = 0$ and
$\infty$, there is one singular fiber of type $\IV$ at $u_1=1$. It is
the divisor $E_{3,2}+E'_{3,2}+Q_1$ (the long dashed dotted lines in
Figure~\ref{no1}), where $Q_1$ is a $(-2)$-curve on $X_3$ arising from a
curve on $\mathbb{P}^1 \times \mathbb{P}^1$ below.

Let $p_j : C_{\omega}^j \to \mathbb{P}^1 \; (j=1,2)$ be the projection
given by
\begin{equation}
  \begin{array}{cccc}\label{proj}
    p_j : & C_{\omega}^j & \to & \mathbb{P}^1\\
    & (x_j :  y_j : z_j) & \mapsto &
     \begin{cases}
       (y_j : z_j) \quad \text{if $z_j \neq 0$}\\
       (1:0) \quad \text{if $z_j = 0$}.
     \end{cases}
  \end{array}
\end{equation}
Then the map $p_1 \times p_2 : C_{\omega}^1 \times C_{\omega}^2 \to
\mathbb{P}^1 \times \mathbb{P}^1$ factors through $\bar{\pi} :
C_{\omega}^1 \times C_{\omega}^2 / \sigma \to \mathbb{P}^1 \times
\mathbb{P}^1$. Let $\pi$ be the morphism of degree three from $X_3$ to
$\mathbb{P}^1 \times \mathbb{P}^1$ that makes the following diagram commutative:
\[\xymatrix{
& X_3 \ar[d]  \ar[dr]^{\pi} & \\
C_{\omega}^1 \times C_{\omega}^2 \ar[r] & C_{\omega}^1 \times C_{\omega}^2 / \sigma \ar[r]_{\bar{\pi}} & \mathbb{P}^1 \times \mathbb{P}^1
}\]
It is easy to verify that the equation $u_1=1$ means
\begin{equation}\label{Q_1}
  y_1^2-2y_1-2y_2-1=0
\end{equation}
from (\ref{par1}).  This equation defines a curve on $\mathbb{P}^1
\times \mathbb{P}^1$. Then it lifts to the $(-2)$-curve $Q_1$ on $X_3$
via the map $\pi$.

\section{Fibration 3}\label{fib3}

\begin{figure}[h]
  \centering
  \includegraphics[width=10cm]{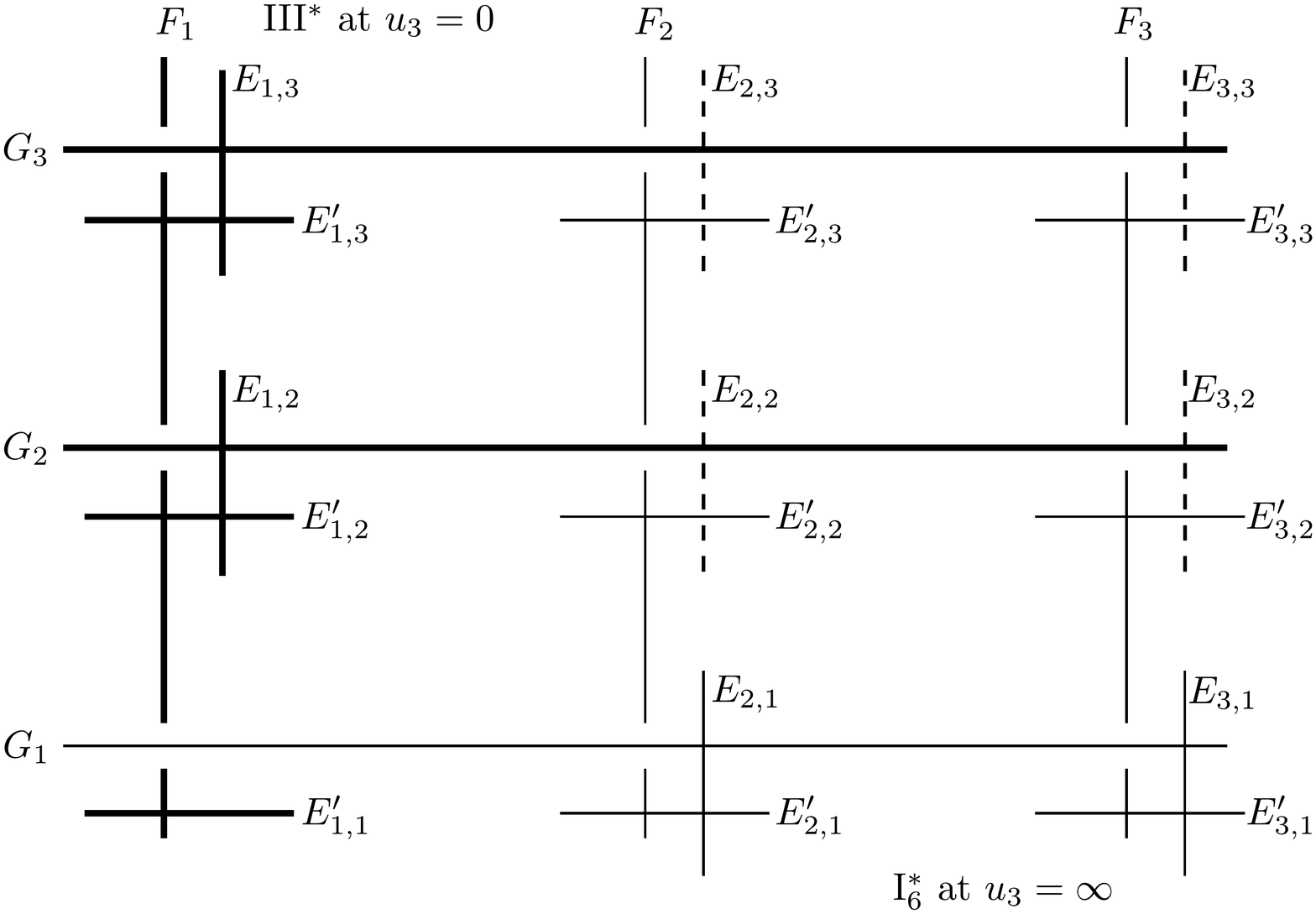}
  \caption{Fibration 3}
  \label{no3}
\end{figure}

An elliptic parameter for Fibration 3 is given by
\begin{equation}
  \label{par3}
  u_3=\frac{t}{y_1^2-1}.
\end{equation}
The divisor of $u_3$ is given by
\begin{equation}
  \label{div3}
  \begin{aligned}
    \left( u_3 \right) &= G_2+2E_{1,2}+3E'_{1,2}+4F_1+3E'_{1,1}+2E_{1,3}+G_3+3E'_{1,2}\\
    & - \left( E'_{2,2} +E'_{2,3}+2 \left( F_2 + E'_{2,1}+E_{2,1}+G_1+E_{3,1}+E'_{3,1}+F_3 \right) + E'_{3,2}+E'_{3,3} \right),
  \end{aligned}
\end{equation}
which is indicated in Figure~\ref{no3}. The zero divisor $(u_3)_0$ is
the singular fiber of type $\III^*$ (the bold lines) and the polar
divisor $(u_3)_{\infty}$ is the singular fiber of type $\I_6^*$ (the
thin lines). The curves $E_{2,2}, E_{2,3}, E_{3,2}$ and
$E_{3,3}$ (the dotted lines) are all the sections.

Eliminating the variable $t$ from (\ref{model}) and (\ref{par3}), we have the following equation
\begin{equation}
  \label{eq3-1}
  y_2^2=u_3^3(y_1^2-1)^2 +1,
\end{equation}
which has a rational point $(y_1, y_2)=(1,1)$ corresponding to the
curve $E_{2,2}$.  Thus, choosing $E_{2,2}$ as the zero section of the
group structure, we obtain the Weierstrass equation for Fibration 3
\begin{equation}
  \label{eq-no3}
  Y^2 = X^3 + 4 u_3^3 X^2 - 4 u_3^3 X,
\end{equation}
where the change of variables is given by
\begin{equation}
  \label{change3}
  X=\frac{2(y_2+1)}{(y_1-1)^2}, \quad Y=\frac{4\left( u_3^3 (y_1+1)(y_1-1)^2+y_2+1 \right)}{(y_1-1)^3}.
\end{equation}
Besides the above two singular fibers of types $\III^*$ and $\I_6^*$, the
fibration has three $\I_1$ fibers at $u_3= -1, -\omega$ and $ -\omega^2$.

The 2-torsion rational point $(X,Y)=(0,0)$ corresponds to the curve
$E_{3,3}$. The rational point $(X,Y)=(1,-1)$ corresponds to the curve
$E_{3,2}$ of height $\langle E_{3,2}, E_{3,2} \rangle =\frac{3}{2}$,
which is a generator of the Mordell-Weil lattice of the fibration. The
curve $E_{2,3}$ is another free section corresponding to the rational
point $(1,1)$ with the relation $E_{2,3}=-E_{3,2}$ in the Mordell-Weil
group.

\section{Fibration 5}\label{fib5}

\begin{figure}[h]
  \centering
  \includegraphics[width=10cm]{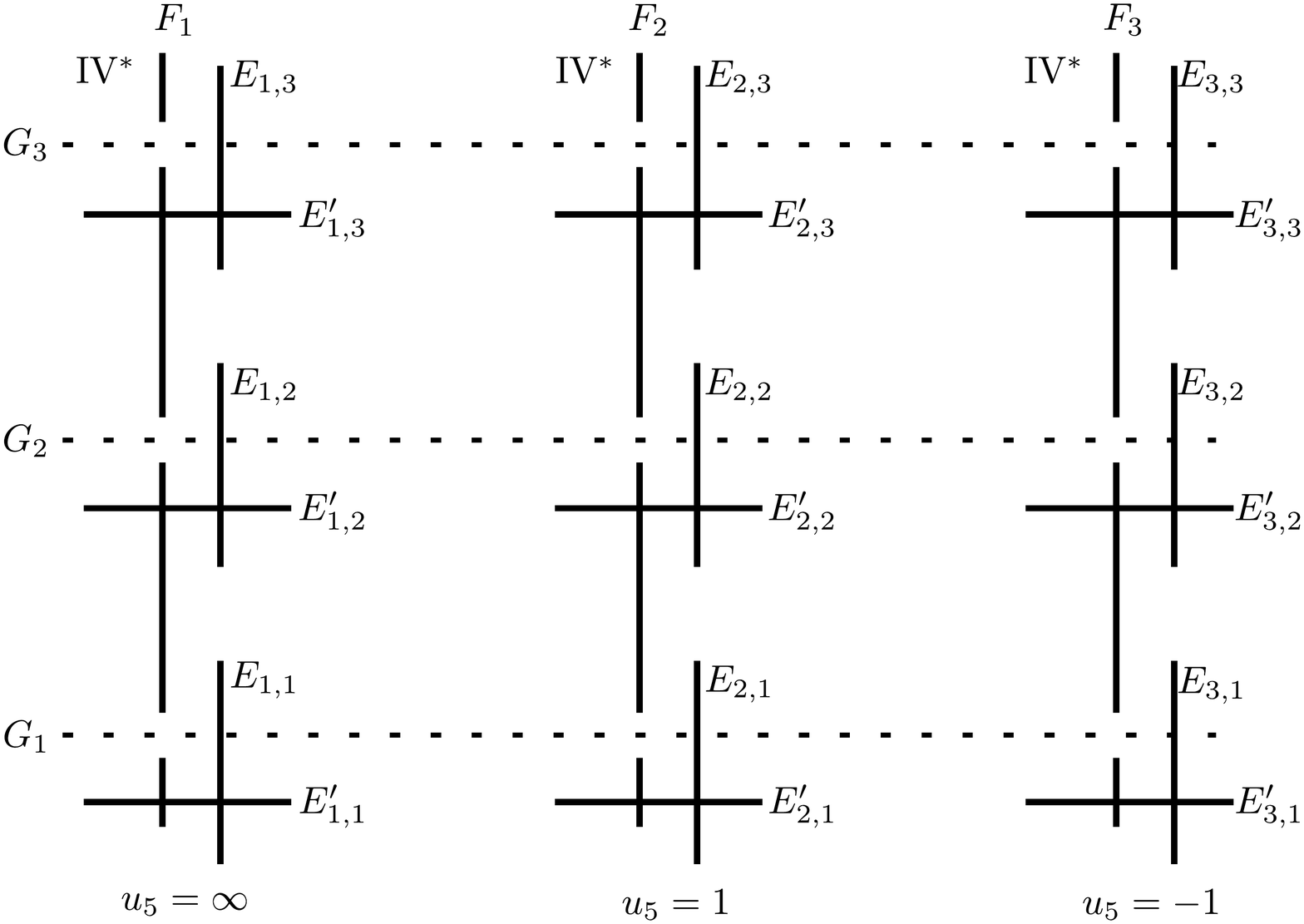}
  \caption{Fibration 5}
  \label{no5}
\end{figure}

An elliptic parameter for Fibration 5 is given by 
\begin{equation}
  \label{par5}
  u_5=y_1.
\end{equation}
It is clear that this elliptic parameter defines a fibration having
three singular fibers all of types $\IV^*$ at $u_5=1, -1$ and $\infty$
(the bold lines in Figure~\ref{no5}) from (\ref{func}). Furthermore
the fibration is induced by the composition of the first projection
$C^1_{\omega} \times C^2_{\omega} \to C^1_{\omega}$ and the covering
map of degree three $p_1: C^1_{\omega} \to \mathbb{P}^1$ in
(\ref{proj}).

The following simple coordinate change
\begin{equation}
  \label{change5}
  X=(u_5^2-1) t, \quad Y= (u_5^2-1)^2 y_2
\end{equation}
 converts the equation (\ref{model}) into the Weierstrass equation for Fibration 5
 \begin{equation}
   \label{eq-no5}
   Y^2= X^3 + (u_5^2-1)^4.
 \end{equation}
 The curve $G_1, \, G_2$ and $G_3$ correspond to the zero section,
 3-torsion rational points $(0, (u_5^2-1)^2)$ and $(0, -(u_5^2-1)^2)$,
 respectively (the dotted lines in Figure~\ref{no5}).

\section{Fibration 6}\label{fib6}

\begin{figure}[h]
  \centering
  \includegraphics[width=10cm]{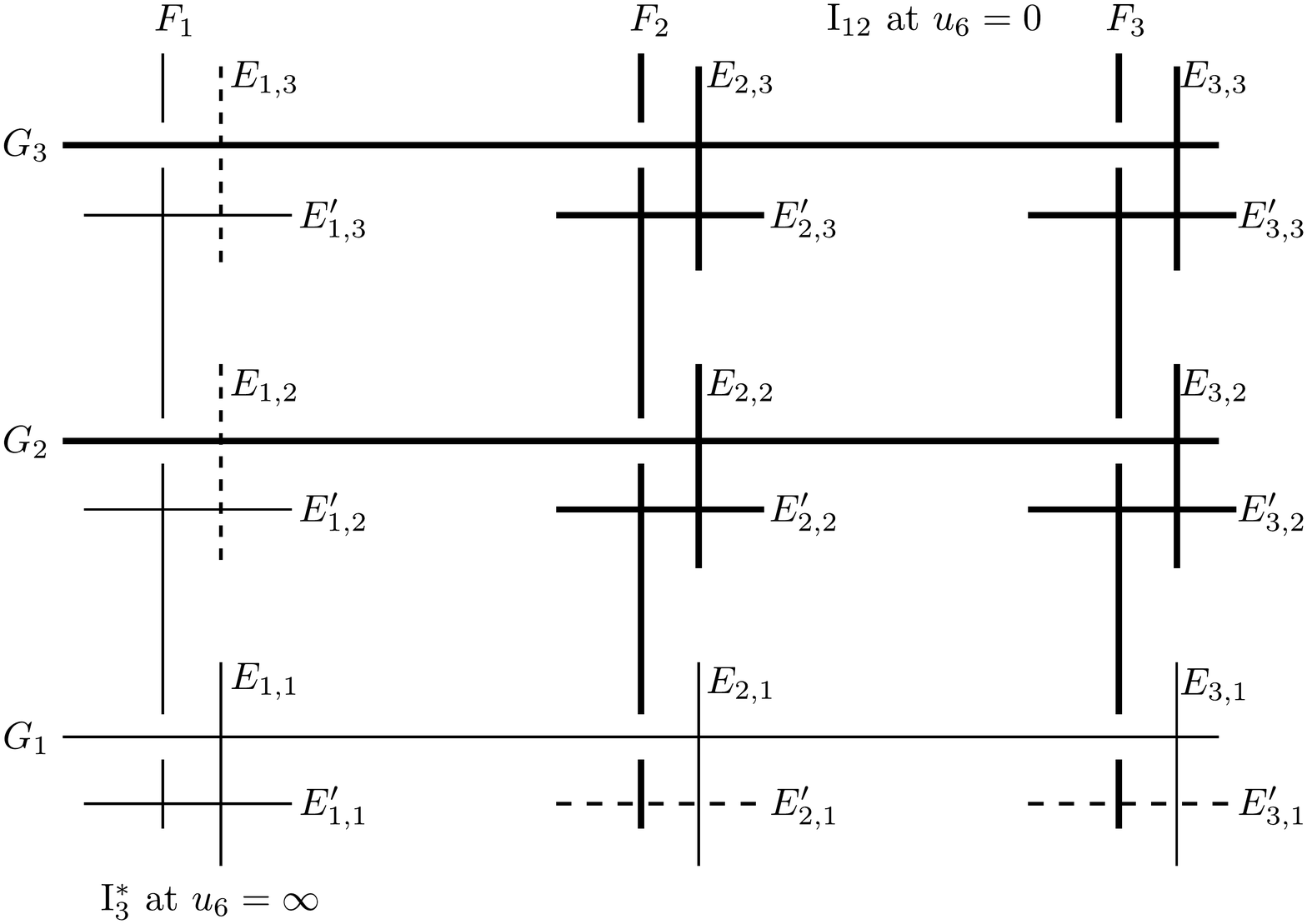}
  \caption{Fibration 6}
  \label{no6}
\end{figure}

An elliptic parameter for Fibration 6 is given by
\begin{equation}
  \label{par6}
  u_6=t.
\end{equation}
Since we gave the divisor of $t$ in (\ref{func}), we know that the
zero divisor $(u_6)_0$ is the singular fiber of type $\I_{12}$ (the
bold lines in Figure~\ref{no6}) and the polar divisor $(u_6)_{\infty}$
is the singular fiber of $\I_3^*$ (the thin lines in
Figure~\ref{no6}). The curves $E_{1,2}, E_{1,3}, E'_{2,1}$ and
$E'_{3,1}$ (the dotted lines in Figure~\ref{no6}) are all the
sections. Choosing $E_{1,2}$ as the zero section of the group
structure, we obtain the Weierstrass equation for Fibration 6
\begin{equation}
  \label{eq-no6}
  Y^2 = X^3 -2(u_6^3-2)X^2-u_6^6X,
\end{equation}
where the change of variables is given by
\begin{equation}
  \label{change6}
  X=\frac{t^3(y_2+1)}{y_2-1}, \quad Y=\frac{2t^3y_1(y_2+1)}{y_2-1}.
\end{equation}
Besides the two singular fibers of type $\I_{12}$ at $u_6=0$ and of
type $\I_3^*$ at $u_6 = \infty$, there are three $\I_1$ fibers at
$u_6=1, \omega$ and $\omega^2$. The Mordell-Weil group of the
fibration is isomorphic to $\mathbb{Z}/4\mathbb{Z}$. The curve
$E_{1,3}$ corresponds to the rational point $(0,0)$ of order two, and
remaining curves $E'_{2,1}$ and $E'_{3,1}$ correspond to the rational
points $(u_6^3,2u_6^3), \, (u_6^3,-2u_6^3)$ of order four,
respectively.

\section{Fibration 4}\label{fib4}

To obtain the Weierstrass equation for Fibration 4, we use a
2-neighbor step from Fibration 3. For more detail about
\textit{2-neighbor step}, we refer to \cite{Kumar,Sengupta, Utsumi}.
\begin{figure}[h]
  \centering
  \includegraphics[width=10cm]{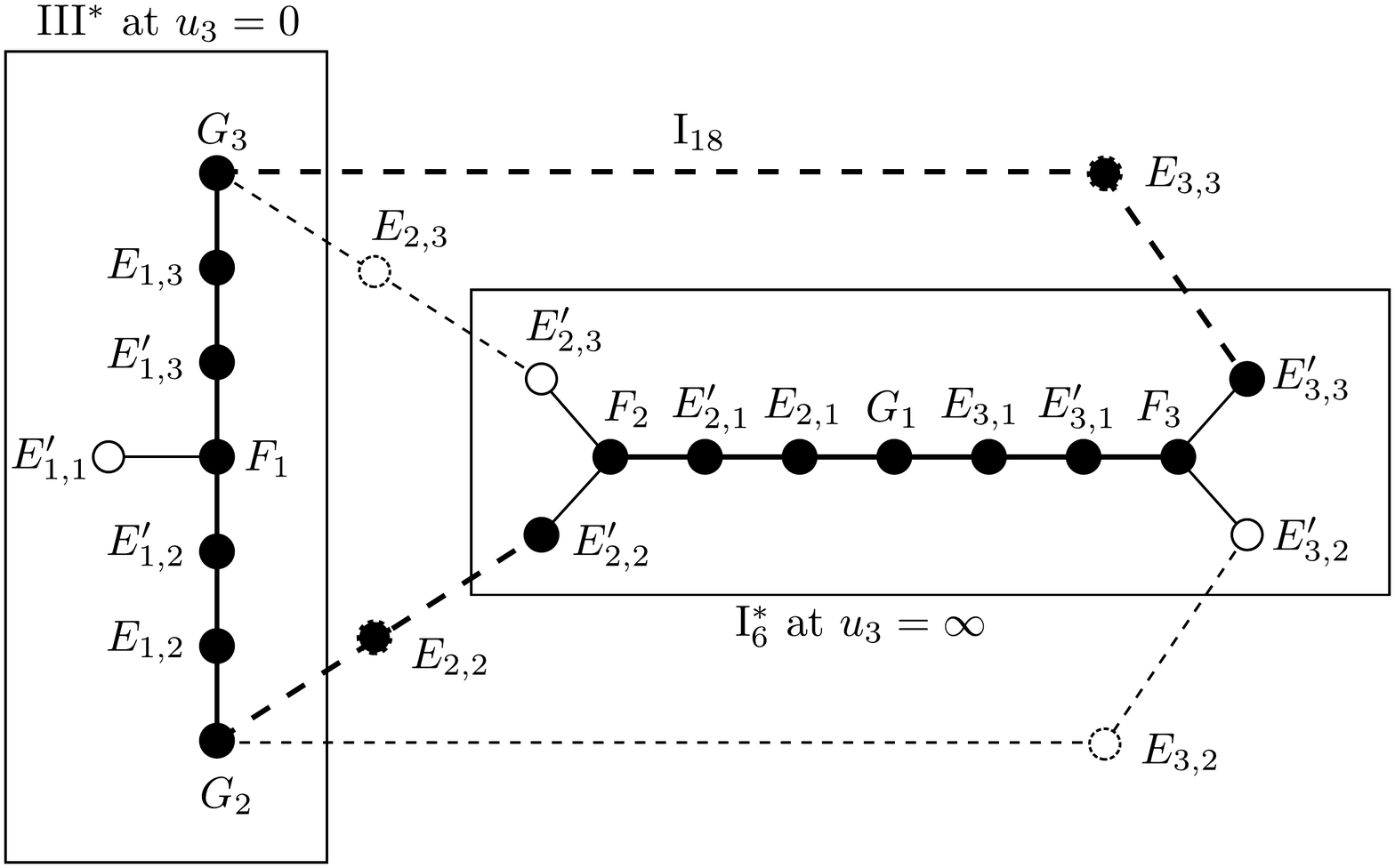}
  \caption{2-neighbor from Fibration 3 to Fibration 4}
  \label{no3-4}
\end{figure}

We compute explicitly the elements of $\mathcal{O}_{X_3}(F)$ where
\begin{equation}\label{div-no4}
  \begin{aligned}
  F &= E_{2,2}+ G_2 + E_{1,2} + E_{1,2}' + F_1 + E_{1,3}' + E_{1,3} + G_3 + E_{3,3} + E_{3,3}' + F_3\\
  & \quad + E_{3,1}' + E_{3,1} + G_1 + E_{2,1} + E_{2,1}' + F_2 + E_{2,2}'
\end{aligned}
\end{equation}
is the class of the fiber of type $\I_{18}$ we are considering. The
linear space $\mathcal{O}_{X_3}(F)$ is 2-dimensional, and the ratio of
two linearly independent elements is an elliptic parameter for
$X_3$. Since $1$ is an element of $\mathcal{O}_{X_3}$, we may find a
non-constant element of $\mathcal{O}_{X_3}(F)$. Then it will be an
elliptic parameter of Fibration 4. Let us $u_4' \in
\mathcal{O}_{X_3}(F)$ be a non-constant. The function $u_4'$ has a
simple pole along $E_{2,2}$ and $E_{3,3}$, which are the zero section
and 2-torsion of Fibration 3. Also, it has a simple pole along $G_2$,
the identity component of the fiber at $u_3=0$, a simple pole along
$E_{3,3}'$, the identity component of the fiber at
$u_3=\infty$. Therefore we can put
\begin{equation}
  u_4' = \frac{\frac{Y}{X} + A_0 + A_1 u_3 + A_2 u_3^2}{u_3},
\end{equation}
where the variables $u_3, X, Y$ are given by (\ref{par3}) and
(\ref{change3}). Assume $A_1=0$, since $1$ is an element of
$\mathcal{O}_{X_3}(F)$. To obtain the coefficients $A_0$ and $A_2$, we
look at the order of vanishing along the non-identity components of
fibers at $u_3=\infty$. The function $u_4'$ does not have any pole
along $E_{3,2}'$, which intersects with the section $E_{3,2}$ of the
fibration 3 at $u_3 = \infty$. Hence $u_4'$ has no pole at $(X,Y,u_3)
= (1,-1,\infty)$, and that gives us $A_2=0$. Similarly, the component
$E_{2,3}'$, which intersects with the section $E_{2,3}$, gives us
$A_0=0$. Consequently, we have a new elliptic parameter
\begin{equation}
  \label{par4-neighbor}
  u_4' = \dfrac{Y}{u_3X},
\end{equation}
where the variables $u_3, X, Y$ are given by (\ref{par3}) and
(\ref{change3}).  Solving for $Y$ and substituting into the
Weierstrass equation (\ref{eq-no3}), after suitable coordinate changes
we have the following
\begin{equation}
  \label{eq-no4-0}
  y^2 = x^3+\frac{1}{4}({u_4'}^2 x-16)^2. 
\end{equation}
Although this is a Weierstrass equation for Fibration 4, for latter
calculations, we put
\begin{equation}
  \label{basechange4}
  {u_4}'=\frac{2}{u_4}, \; x=\frac{2^2X}{u_4^4}, \; y=\frac{2^3Y}{u_4^6}
\end{equation}
and obtain another Weierstrass equation for Fibration 4
\begin{equation}
  \label{eq-no4}
  Y^2=X^3+(X-u_4^6)^2.
\end{equation}
The change of variables is given by 
\begin{equation}
  \label{change4}
  u_4 = \dfrac{t}{y_1+y_2}, \; X=\frac{(y_1^2-1)t^3}{(y_1+y_2)^4}, \; Y=\frac{(y_1^2y_2+2y_1+y_2) t^6}{(y_2^2-1)(y_1+y_2)^6}.
\end{equation}
The fibration has singular fibers of type $\I_{18}$ at $u_4=0$ and of
type $\I_1$ at the zeros of $27u_4^6+4=0$. The zero section
corresponds to the divisor $E'_{1,1}$. The 3-torsion rational points
$(0,u_4^6)$ and $(0,-u_4^6)$ correspond to the divisors $E'_{3,2}$ and
$E'_{2,3}$, respectively. The free rational points $(2u_4^3,
u_4^4+2u_4^3)$ and $(-2u_4^6,u_4^3-2u_4^3)$ correspond to the divisors
$E_{3,2}$ and $E_{2,3}$, respectively with the relation $E_{2,3} +
E_{3,2} = E'_{2,3}$ in the Mordell-Weil group. Since the height of
$E_{2,3}$ is equal to $\frac{3}{2}$, $E_{2,3}$ generates the
Mordell-Weil lattice of the fibration.

\begin{figure}[h]
  \centering
  \includegraphics[width=10cm]{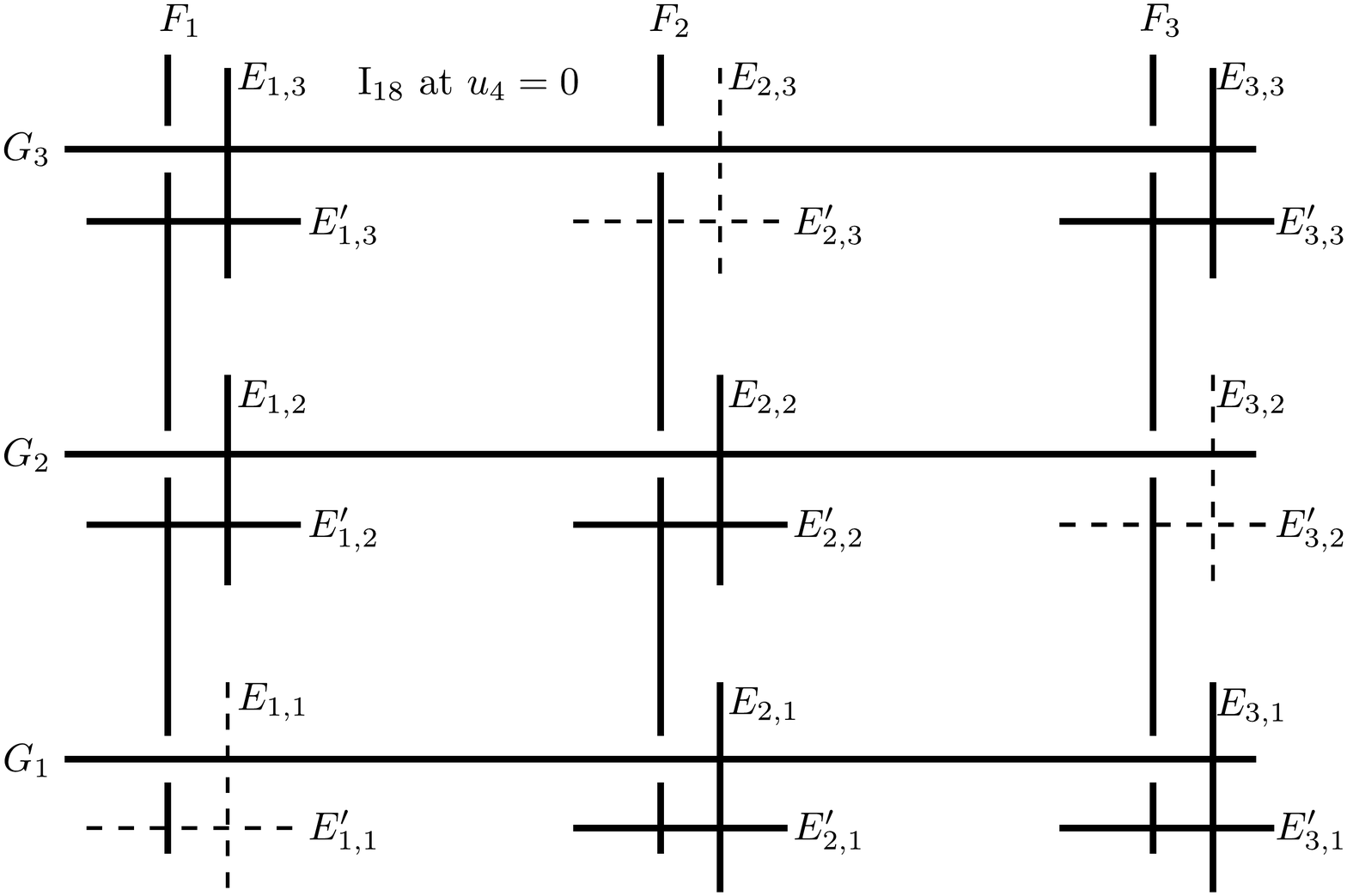}
  \caption{Fibration 4}
  \label{no4}
\end{figure}

\section{Fibration 2}\label{fib2}

\begin{figure}[h]
  \centering
  \includegraphics[width=10cm]{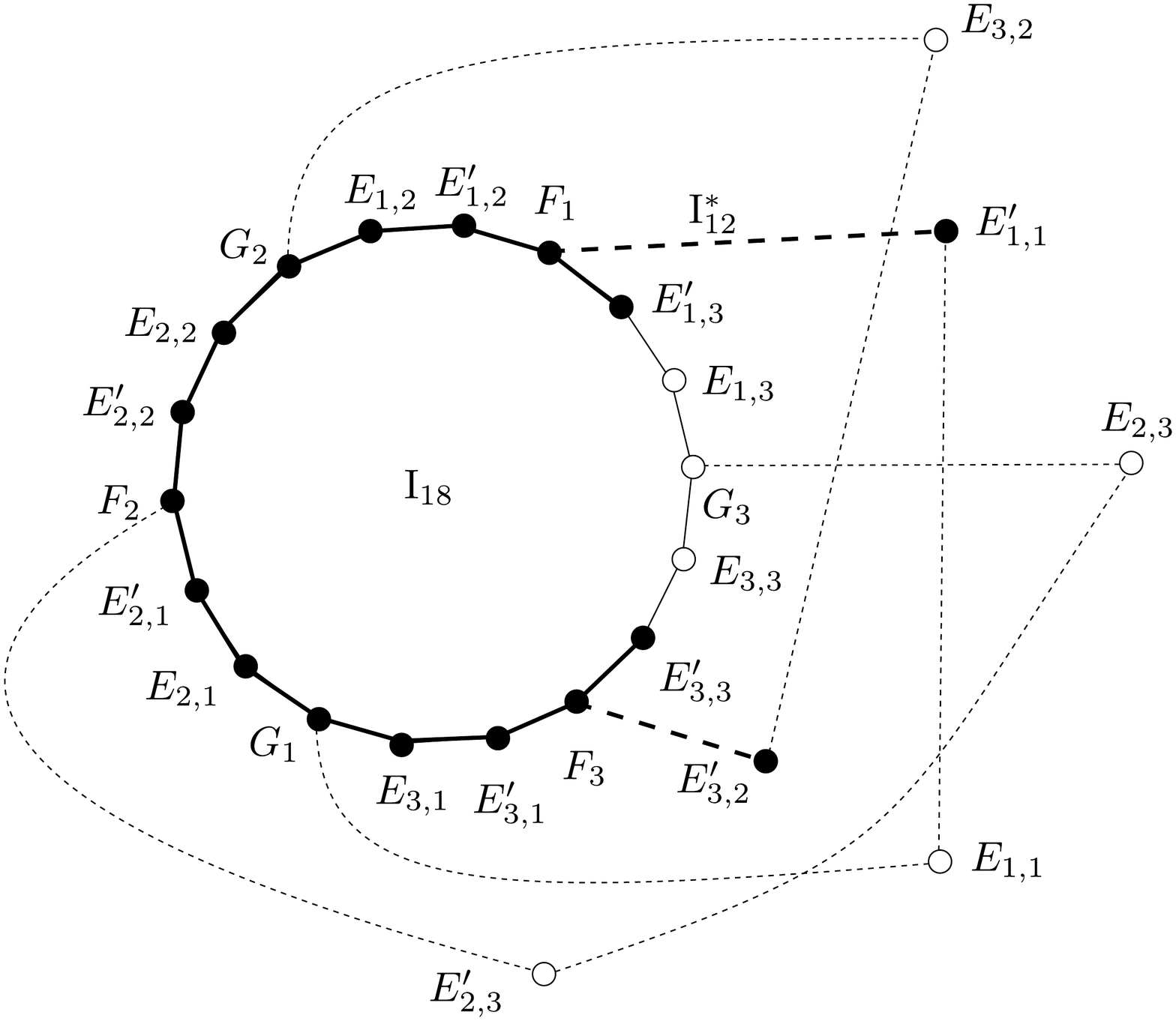}
  \caption{2-neighbor form Fibration 4 to Fibration 2}
  \label{no4-2}
\end{figure}

\begin{figure}[h]
  \centering
  \includegraphics[width=10cm]{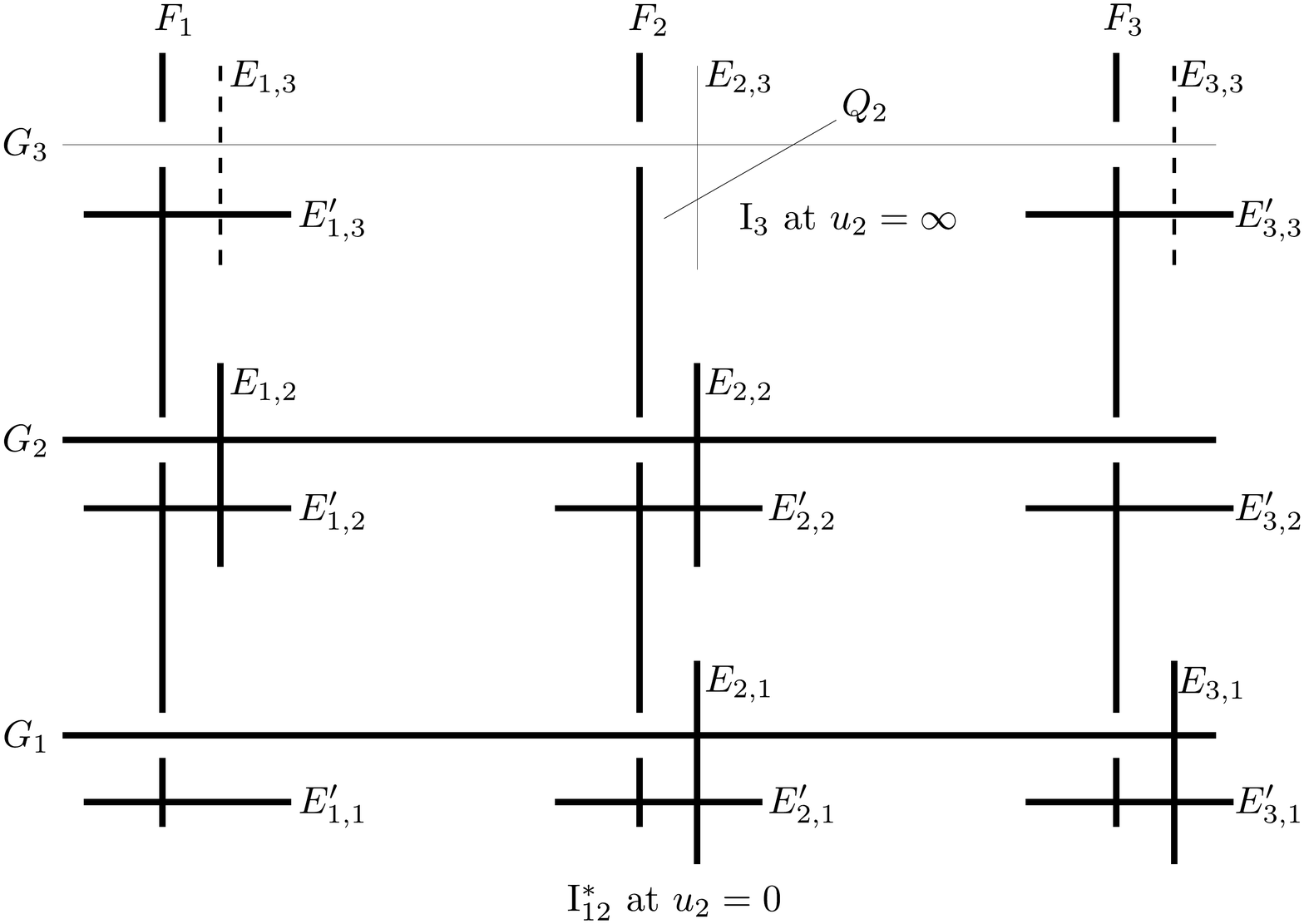}
  \caption{Fibration 2}
  \label{no2}
\end{figure}

We obtain the following elliptic parameter $u_2'$ for Fibration 2 by a
2-neighbor step from Fibration 4 (see Figure~\ref{no4-2}).
\begin{equation}
  \label{par2-0}
  u_2' = \frac{u_4^6 + X + Y}{u_4^2 X}
\end{equation}
The variables $u_4, X, Y$ are given by (\ref{change4}). 
Then we get the following Weierstrass equation for Fibration 2.
\begin{equation}
  \label{eq-no2-0}
  y^2=x^3+2({u_2'}^3-4)x^2+16x.
\end{equation}
We put
\begin{equation}
  \label{basechange2}
  u_2'=\frac{2}{u_2}, \; x=\frac{2^2 X}{u_2^4}, \; y=\frac{2^3 Y}{u_2^6}
\end{equation}
and obtain another Weierstrass equation for Fibration 4.
\begin{equation}
  \label{eq-no2}
  Y^2=X^3-2(u_2^3-2)X^2-u_2^8X.
\end{equation}
The change of variables is given by
\begin{equation}
  \label{change2}
  \begin{aligned}
    &u_2=\frac{2t^2}{(y_2+1)(y_1^2+2y_1+2y_2-1)},\\
    &X=-\frac{32(y_1-1)^2(y_2-1)^3t^2}{(y_2+1)^2(y_1^2+2y_1+2y_2-1)^4},\\
    &Y= -\frac{128(y_1-1)^3(y_2-1)^4(y_1+1)(y_1+y_2)}{(y_2+1)^2(y_1^2+2y_1+2y_2-1)^5}.
  \end{aligned}
\end{equation}
The zero divisor $\left( u_4 \right)_0$ is the singular fiber of type
$\I_{12}^*$ (the bold lines in Figure~\ref{no2}). The polar divisor
$\left(u_4\right)_{\infty} = G_3 + E_{2,3} + Q_2$ is the singular
fiber of type $\I_3$ (the thin lines in Figure~\ref{no2}), where the
divisor $Q_2$ is the lifting of the curve $y_1^2+2y_1+2y_2-1=0$ on
$\mathbb{P}^1 \times \mathbb{P}^1$ by the map $\pi$ in \S\ref{fib1}.
Besides these two singular fibers, there are three $\I_1$ fibers at
$u_2=1, \omega$ and $\omega^2$. The zero section corresponds to the
divisor $E_{1,3}$. The 2-torsion rational point $(0,0)$ corresponds to
the divisor $E_{3,3}$.

\begin{remark}
  We give a Weierstrass equation for Fibration 6 in
  \S\ref{fib6}. Comparing the equations (\ref{eq-no2}) and
  (\ref{eq-no6}), we know easily that Fibration 2 is a quadratic
  twist of Fibration 6. This is the reason why we adopt the equation
  (\ref{eq-no2}) as the Weierstrass equation for Fibration 2 rather
  than the equation (\ref{eq-no2-0}).
\end{remark}

\begin{acknowledgements}
  The computer algebra system Maple and Maple Library ``Elliptic
  Surface Calculator'' written by Professor Masato Kuwata~\cite{ESC}
  were used in the calculation for this paper. The author would like
  to thank the developers of these programs.
\end{acknowledgements}

\end{document}